\documentclass{article}
\usepackage{amssymb,graphicx}
\newtheorem{theorem}{Theorem}[section]

\newtheorem{lemma}[theorem]{Lemma}

\newtheorem{corollary}[theorem]{Corollary}
\def\binom#1#2{{#1}\choose{#2}}

\newcommand{\dd}{{, \ldots ,}}
\newcommand{\df}{\displaystyle\frac}

\newcommand{\la}{\lambda}

\newcommand{\ga}{\gamma}

\newcommand{\qed}{\square}

\newcommand{\RR}{\mathbb{R}}
\newcommand{\ZZ}{\mathbb{Z}}
\newcommand{\bw}{{\bf w}}
\newcommand{\bv}{{\bf v}}
\newcommand{\be}{{\bf e}}
\newcommand{\bn}{{\bf n}}
\newcommand{\bz}{{\bf z}}

\newcommand{\eps}{\epsilon}

\newcommand{\bS}{{\bf{S}}}
\newcommand{\bM}{{\bf{M}}}
\newcommand{\bN}{{\bf{N}}}

\newcommand{\sL}{{\mathcal L}}
\newcommand{\sC}{{\mathcal C}}
\newcommand{\sD}{{\mathcal D}}
\newcommand{\sG}{{\mathcal G}}
\newcommand{\sE}{{\mathcal E}}
\newcommand{\sH}{{\mathcal H}}
\newcommand{\sK}{{\mathcal K}}

\newcommand{\sT}{{\mathcal T}}
\newcommand{\sW}{{\mathcal W}}

\usepackage{amssymb,graphicx}

\begin{document}
\title{The Computational Complexity of Knot and
Link Problems}
\author{ Joel Hass \footnote {Supported in part by NSF grant DMS-9704286.},
Jeffrey C. Lagarias
\footnote{ Work done in part while the first two authors were visiting the
Mathematical Sciences Research Institute in Berkeley 1996/7. Research at
MSRI is supported in part by NSF grant DMS-9022140.}
and Nicholas Pippenger
\footnote{Supported in part by NSERC research grant OGP 0041640.}.
}
\maketitle

\begin{abstract}
We consider the problem of deciding whether a polygonal knot in 3-dimensional Euclidean space is unknotted, capable of 
being continuously deformed without self-intersection so that it lies in a 
plane. We show that this problem, {\sc unknotting problem} is in {\bf NP}.
We also 
consider the problem, {\sc unknotting problem} of determining whether two or 
more such polygons can be split, or
continuously deformed without self-intersection so that they occupy both 
sides of a plane without intersecting it. We show that it also is in {\bf NP}.
Finally, we show that the problem of determining the genus of a polygonal 
knot (a generalization of the problem of determining whether it is 
unknotted) is in {\bf PSPACE}. We also give exponential worst-case running 
time bounds for deterministic algorithms to solve each of these problems.
These algorithms are based on the use of normal surfaces and decision 
procedures due to W. Haken,
with recent extensions by W. Jaco and J. L. Tollefson. 
\end{abstract}

\vspace*{5\baselineskip}

\newpage 
\section{Introduction}\label{intro}
The problems dealt with in this paper might reasonably be called 
``computational topology''; that is, we study classical problems of 
topology (specifically, the topology of $1$-dimensional curves in
$3$-dimensional space) with the objective of determining their computational 
complexity. One of the oldest and most fundamental of such problems is 
that of determining whether a closed curve embedded in space is 
unknotted (that is, whether it is capable of being continuously deformed 
without self-intersection so that it lies in a plane). Topologists study 
this problem at several levels, with varying meanings given to the terms 
``embedded'' and ``deformed''. The level that seems most appropriate for 
studying computational questions is that which topologists call 
{\em piecewise-linear}. At this level, a closed curve is embedded in space as a 
simple (non-self-intersecting) polygon with finitely many edges.
Such an embedding is called a {\em knot}.
(Operating at the piecewise-linear level excludes {\em wild} knots, such as 
those given by certain polygons with infinitely many edges which
do not have nice neighborhoods or thickenings.)
More generally, one may study {\em links.} A link is a finite collection 
of simple polygons disjointly embedded in 3-dimensional space. The 
individual polygons are called components of the link and a knot is a link 
with one component.

A continuous deformation is required to be piecewise-linear; that is, it 
consists of a finite number of stages, during each of which every vertex of 
the polygon moves linearly with time. 
From stage to stage the number of edges in the polygon may increase 
(by subdivision of edges
at the beginning of a stage) or decrease (when cyclically consecutive 
edges become collinear at the end of a stage). If the polygon remains 
%spelling collinear?
simple throughout this process, the deformation is called an {\em isotopy} 
between the initial and final knots.
Knot isotopy defines an equivalence relation, called {\em equivalence} of 
knots. It is easy to see that all knots that lie in a single plane are 
equivalent; knots in this equivalence class are said to be {\em unknotted} or 
{\em trivial} knots.

We may remark at this point that while it is ``intuitively clear'' that 
there are non-trivial knots, it is not at all obvious how to prove this. 
Stillwell \cite{Stil:80} traces the mathematical notion of knot back to a 
paper of A.~T. Vandermonde in 1771; the first convincing proof of the 
non-triviality of a knot seems to be due to Max Dehn \cite{Dehn:10} in 
1910.

There are a great many alternative 
formulations of the notion of knot equivalence. Here are some.

\begin{enumerate}
\item One can consider sequences of {\em elementary moves,} which are very 
simple isotopies that move a single edge across a triangle to the opposite 
two sides, or vice versa.

\item One can consider {\em ambient isotopies} that move not only the knot, 
but also the space in which it is embedded, in a piecewise-linear way. 

\item One can consider {\em homeomorphisms} (continuous bijections that 
have continuous inverses) that map the space to itself in a piecewise-linear way,
are orientation preserving, and send one knot to the other. 
\end{enumerate}

\noindent%
One can also study knots or links by looking at their {\em projections} onto a 
generic plane.
In this way, a knot or link may be represented by a planar graph, called a 
{\em knot diagram} or {\em link diagram}, in which all vertices (representing the 
{\em crossings} of edges of the polygon) have degree four, and for which an 
indication is given at each crossing of which edge goes ``over'' and which 
edge goes ``under''. This gives an additional formulation of equivalence: 

\begin{enumerate}
\setcounter{enumi}{3}
\item One may consider
sequences of {\em Reidemeister moves}, which are simple transformations on 
the diagram of a knot that leave the equivalence class of the knot 
unchanged.
\end{enumerate}

\noindent%
For more details on piecewise-linear topology, the various formulations 
of knot and link equivalence, and many other aspects of knot theory, we 
recommend the books \cite{Adams}, \cite{Bur:85}, \cite{Rolfsen}. 
An introduction to the notions of complexity which we use can be found in
\cite{Garey-Johnson},\cite{Wel:93a}.

In order to study the computational complexity of knot and link problems, 
we must agree on a finite computational representation of a knot or link. 
There are two natural representations: a polygonal representation in 3-dimensional space,
or a link diagram representing a 2-dimensional projection. 

A polygonal representation of a link $L$ consists of a set of simple 
polygons in 3-dimensional space described by listing the vertices of each 
polygon in order; we assume that these vertices have rational coordinates. 
We can reduce to the case of integer lattice point vertices by replacing 
$L$ by a scaled multiple $mL$ for a suitable integer $m$. This does not 
change the equivalence class of $L$. A particularly simple kind of 
polygonal representation uses only integer lattice points as vertices and 
edges of unit length, so that the polygon is a closed self-avoiding walk on 
the integer lattice; a sequence of moves (up, down, north, south, east, 
west) that traverse the polygon, returning to the starting point without 
visiting any other point twice. (This formulation was used by Pippenger 
\cite{Pip:89} and Sumners and Whittington \cite{SW:88} to show that 
``almost all''
long self-avoiding polygons are non-trivially knotted.) The size of a 
polygonal representation $L$ is the number of edges in $L$; its input 
length is the number of bits needed to describe its vertices, in binary.

A link diagram ${\cal D}$ is a planar graph with some extra labeling for 
crossings that specifies a (general position) two-dimensional projection 
of a link.
A precise definition is given in
Section 3. The size of a link diagram is the number of vertices in $D$.

These two representations are polynomial-time equivalent in the 
following sense.
Given a polygonal representation $L$ one can find in polynomial time in its 
input length a planar projection yielding a link diagram ${\cal D}$; if $L$ 
has $n$ edges then the graph ${\cal D}$ has at most
$O(n^2)$ vertices. Conversely given a link diagram ${\cal D}$ with $n$ 
vertices
and $l$ components, one can compute in time polynomial in $n+l$ a 
polygonal link $L$ with $O(n+l)$ edges that has integer vertices and input 
length $O(n+l)$ and which projects in the $z$-direction onto the link 
diagram ${\cal D}$; see Section 7.

In this paper we consider knots and links as represented by link diagrams 
and take the crossing number as the measure of input size. We can now 
formulate the computational problem of recognizing unknotted polygons as 
follows: 

\begin{tabbing}
{\em Problem:} UNKNOTTING PROBLEM \\
{\em Instance:} A link diagram ${\cal D}$. \\
{\em Question:} \= Is ${\cal D}$ a knot diagram that represents the trivial knot? \\
\end{tabbing}

See Welsh \cite{Wel:93a}---\cite{Wel:93c} for more information on this 
problem.
The main result of this paper is the following. 

\begin{theorem}
\label{knotnp}
The {\sc unknotting problem} is in {\bf NP}. \end{theorem} 

The {\sc unknotting problem} was shown to be decidable by Haken 
\cite{Hak:61}; the result was announced in 1954, and the proof published 
in 1961. From then until now, we know of no strengthening of Haken's 
decision procedure to give an explicit complexity bound.
We present such a bound in Theorem~\ref{Nth81}.

We also study the splittability of links. A link is said to be {\em splittable} if 
it can be continuously deformed (by a piecewise-linear
isotopy so that one or more curves of the link can be separated from one 
or more other curves by a plane that does not itself intersect any of the 
curves.
We note that this notion remains unchanged if we replace ``plane'' by 
``sphere'' in the definition.
We formulate the computational problem of recognizing splittable links as 
follows.

\begin{tabbing}
{\em Problem:} SPLITTING PROBLEM \\
{\em Instance:} A link diagram ${\cal D}$. \\ {\em Question:} Is the link 
represented by ${\cal D}$ splittable? \end{tabbing} 

The {\sc splitting problem} was shown to be decidable by Haken \cite{Hak:61} 
in 1961, see also Schubert \cite{Sch:61}. We establish the following 
result. 

\begin{theorem}
\label{linknp}
The {\sc splitting problem} is in {\bf NP}. \end{theorem} 

Another generalization of the unknotting problem concerns
the {\em genus} $g(K)$ of a knot $K$.
This is an integer assoicated to each knot, which
is invariant under isotopy. It was defined 
by Seifert \cite{Sei:35} in 1935; an informal account of the definition
follows.
Given a knot $K$, consider the class ${\cal S}(K)$ of all orientable 
spanning surfaces for $K$; that is, embedded orientable surfaces that have 
$K$ as their boundary. Seifert showed that this class is non-empty for any 
knot $K$. We assume in this discussion that all surfaces are 
triangulated and embedded in a piecewise-linear way. Up to piecewise-linear
homeomorphism, an orientable surface is characterized by the 
number of boundary curves and the number of ``handles'', which is called 
the {\em genus} of the surface.
The genus $g(K)$ of the knot $K$ is defined to 
be the minimum genus of any surface in ${\cal S}(K)$.
Seifert showed that a trivial knot $K$ is characterized by the condition 
$g(K)=0$. Since an orientable
surface with one boundary curve and no handles is homeomorphic to a disk, 
this means that a knot is trivial if and only if it has a spanning disk. 

The notion of genus gives us a natural generalization of the problem of 
recognizing unknotted polygons; we formulate the problem of computing 
the genus as a language-recognition problem in the usual way. 

\begin{tabbing}
{\em Problem:} GENUS PROBLEM \\
{\em Instance:} A link diagram ${\cal D}$ and a natural number $k$. \\ 
{\em Question:} \=Does the link diagram ${\cal D}$ represent a knot $K$ with $g(K) \le k$?
\end{tabbing}

Haken \cite{Hak:61} also gave a decision procedure to determine the genus 
of a knot. We establish the following result.

\begin{theorem}
\label{genusps}
The {\sc genus problem} is in {\bf PSPACE}. \end{theorem} 

In 1961 Schubert \cite{Sch:61} extended Haken's methods to show the 
decidability of the {\sc genus problem} on an arbitrary compact 3-manifold 
with boundary. (In this more general situation, it is necessary to first 
determine whether there exists any orientable surface which has the 
given embedded knot as boundary, and if so, determine the minimal genus.) 
Our analysis in principle extends to this case. 

We also obtain exponential worst case running time bounds for 
deterministic algorithms to solve the above three problems. The input size 
$n$ is measured by the crossing measure of a link diagram. For the 
{\sc unknotting problem} and the {\sc splitting problem} these algorithms run in 
time $O(2^{cn})$ and space $O(n \log^2 n)$ on a Turing machine. We show 
that the genus $g(K)$ for a knot $K$ presented as a knot diagram ${\cal K}$ can 
be computed in time $O(2^{c n^2})$ and space $O(n^2)$, on a Turing 
machine. 

The results in this paper were announced in the Proceedings of the 
38th Annual Symposium on Foundations of Computer Science in Miami Florida in 
October, 1997 \cite{HLP}.

\section{Historical Background}
Recognizing whether two knots are equivalent has been one of the motivating problems 
of knot theory.
A great deal of effort has been devoted to a quest for algorithms for 
recognizing the unknot,
beginning with the work of Dehn \cite{Dehn:10} in 1910. Dehn's idea was 
to look at the fundamental group of the complement of the knot, for which 
a finite presentation in terms of generators and relations can easily be 
obtained from a standard presentation of the knot. Dehn claimed that a 
knot is trivial if and only if the corresponding group is infinite cyclic. The 
proof of what is still known as ``Dehn's Lemma'' had a gap, which remained 
until filled by Papakyriakopoulos in 1957 \cite{Pap:57}. A consequence is 
the criterion that a curve is knotted if and only if the
fundamental group of its complement is nonabelian. Dehn also posed the 
question of deciding whether a finitely presented group is isomorphic to 
the infinite cyclic group. During the 1950's it was shown that many such 
decision problems for finitely presented groups, not necessarily arising from knots, are undecidable (see Rabin 
\cite{Rab:58}, for example),
thus blocking this avenue of progress.
The avenue has been traversed in the reverse direction, however: there 
are decision procedures for restricted classes of finitely presented 
groups arising from topology. In particular, computational results for 
properties of knots that are characterized by properties of the 
corresponding groups can be interpreted as computational results for knot 
groups.

Abstracting somewhat from Dehn's program, we might try to recognize 
knot triviality by finding an invariant of the knot that (1) can be computed 
easily and (2)
assumes some particular value only for the trivial knot. (Here {\em invariant} 
means invariant under isotopy.) Thus Alexander \cite{Alex:28} defined in 
1928 an invariant $A_K(x)$ (a polynomial in the indeterminate $x$) of the 
knot $K$ that can be computed in polynomial time. Unfortunately, it turns 
out that many non-trivial knots have Alexander polynomial $A_K(x) = 1$, 
the same as the
Alexander polynomial of the trivial knot.

Another invariant that has been investigated with the same hope is the 
Jones polynomial $J_K(x)$ of a knot $K$, discovered by Jones \cite{Jon:85} 
in 1985. In this case the complexity bound is less attractive: the Jones 
polynomial for links
(a generalization of the Jones polynomial for knots) is {\bf \#P}-hard and 
in ${\bf FP}^{\bf \#P}$ (see Jaeger, Vertigan and Welsh \cite{JVW:90}). It 
is an open question whether trivial knots are characterized by their Jones 
polynomial. Even this prospect, however, has led Welsh \cite{Wel:93a} to 
observe that an affirmative answer to the last open question would yield 
an algorithm in
${\bf P}^{\bf \#P}$ for recognizing trivial knots, and to add: ``By the 
standards of the existing algorithms, this would be a major advance.'' 

The revolution started by the Jones polynomial has led to the discovery of 
a great number of new knot and link invariants, including Vassiliev 
invariants and invariants associated to topological quantum field
theories, see Birman \cite{Bir:93} and Sawin \cite{Saw:96}. The exact 
ability of these invariants to distinguish knot types has not been 
determined.

A different approach to the problems of recognizing unknottedness and 
deciding knot equivalence eventually led to decision procedures. 
This approach is based on the study of normal surfaces in
3-manifolds (defined in 
Section 3), which was initiated by Kneser \cite{Knes:29}
in 1929. In the 1950's Haken developed
the theory of normal surfaces, and obtained a decision procedure for 
unknottedness in 1961. Haken considered compact orientable 3-manifolds. 
Schubert \cite{Sch:61} extended Haken's procedure to decide the knot 
genus problem and related problems on arbitrary compact 3-manifolds 
with boundary. Haken also outlined an approach via normal surfaces to 
decide the {\sc knot equivalence problem} \cite{Wald:78}:

\begin{tabbing}
{\em Problem:} KNOT EQUIVALENCE PROBLEM \\
{\em Instance}: Two link diagrams ${\sD}_1$ and ${\sD}_2$. \\
{\em Question}: \= Are ${\sD}_1$ and ${\sD}_2$ knot diagrams of equivalent knots? \\
\end{tabbing} 

The final step in this 
program was completed by Hemion \cite{Hem:92} in 1979. This program
actually solves a more general decision problem, concerning a large class 
of 3-manifolds, now called Haken manifolds, which can be cut into 
``simpler'' pieces along certain surfaces (incompressible surfaces), 
eventually resulting in a collection of 3-balls. Knot complements are
examples of Haken manifolds. It gives a procedure
to decide if two Haken manifolds are homeomorphic \cite{JO:84}. Recent 
work of Jaco-Oertel \cite{JO:84} 
and Jaco-Tollefson \cite{JT:95} further simplified some of these 
algorithms.

Apart from these decidability results, there appear to be no explicit 
complexity bounds, either upper or lower, for any of the three problems 
that we study. The work of Haken \cite{Hak:61} and Schubert \cite{Sch:61} 
predates the currently used framework of complexity classes and 
hierarchies. Their algorithms were originally presented in a framework 
(handlebody decompositions) that makes complexity analysis
difficult, but it was recognized at the time that implementation of their 
algorithms
would require at least exponential time in the best case. More recently 
Jaco and others reformulated normal surface theory using piecewise 
linear topology, but did not determine complexity bounds.
Other approaches to knot and link algorithms include methods related to
Thurston's geometrization program for 3-manifolds
(see \cite{Hass:9y} for a survey) and
methods based on encoding 
knots as braids (see Birman-Hirsch \cite{BirHir:97}).
These approaches currently have unknown complexity bounds.

Our results are obtained using a version of normal surface theory as
developed by Jaco-Rubinstein \cite{JacRub89}. Among other things 
we show that Haken's original approach yields an algorithm which 
determines if a knot diagram with $n$ crossings is unknotted in time 
$O(\exp(cn^2))$, and that the improved algorithm of Jaco-Tollefson 
runs in time $O(\exp(cn))$, see Theorem~\ref{Nth81}. The complexity class 
inclusions that we prove require some additional observations. 

\section{Knots and Links}\label{kl}
A {\em knot} is an embedding $f : {\bf S}^1 \to \RR^3$, although it is 
usually identified with its image $K=f({\bf S}^1)$; we are thus considering 
{\em unoriented} knots. A {\em link} with $k$ components is a collection 
of $k$ knots with disjoint images.
An equivalent formulation regards a knot as an embedding in the one-point 
compactification ${\bf S}^3$ of ${\RR^3}$, and we will sometimes use 
this setting.

Two knots $K$ and $K'$ are {\em ambient isotopic} if there exists a 
homotopy $h_t : \RR^3 \to \RR^3$ for $0 \le t \le 1$ such that $h_0$ 
is the identity, each $h_t$ is a homeomorphism, and $h_1(K) = K'$. We 
shall also say in this case that $K$ and $K'$ are {\em equivalent} knots. 
Since we consider {\em unoriented knots}, we will also set two
knots equivalent if they are equal after composing with a reflection of 
$S^1$. Our results work
equally well for oriented or unoriented knots. 

A knot or link is {\em tame} if it is ambient isotopic to a piecewise-linear
knot or link, abbreviated as PL-knot or PL-link,
and also called a polygonal knot or link.
This paper considers PL-knots and links. Given this restriction, we
can without further loss of generality restrict our attention to 
the piecewise-linear settings (see Moise \cite{Moise:52}). 

A {\em regular projection} of a knot or link is an orthogonal projection 
into a plane (say $z=0$) that contains only finitely many multiple points, 
each of which is a double point with transverse crossing. Any regular 
projection of a link gives a {\em link diagram} $\sD$, which is an 
undirected labeled planar graph such that: 

\begin{enumerate}

\item Connected components with no vertices are loops,
simple closed curves disjoint from the rest of the diagram. 

\item Each non-loop edge meets a vertex at each of its two ends (possibly 
the same vertex), and has a label at
each end indicating an overcrossing or undercrossing at that end. 

\item Each vertex has exactly four incident edges, two labeled as 
overcrossings and two labeled as undercrossings, and has a cyclic ordering 
of the incident edges that alternates overcrossings and undercrossings. 

\end{enumerate}
\noindent Conversely, every labeled planar graph satisfying these 
conditions is a link diagram for some link. 

Given a link diagram, if we connect the edges across vertices according to 
the labeling, then the diagram separates into $k$ edge-connected 
components, where $k$ is the number of components in the link. A {\em 
knot diagram} is a link diagram having one component. A {\em trivial} knot 
diagram is a single loop with no vertices. It is straightforward
to see that all knots with the same diagram are isotopic. 

We define the {\em crossing measure}
of a link diagram $\sD$
to be
the number of vertices in the diagram,
plus the number of connected components in the diagram, minus one. 
\begin{equation}
\label{eq201}
n = \# ( \mbox{vertices of} ~ \sL ) + \# (\mbox{connected components}) -1 
~. \end{equation}
For knot diagrams, the crossing measure is equal to the {\em crossing number}, which is the number of vertices in the diagram. A trivial knot 
diagram is the only link diagram with crossing measure zero. All other 
link diagrams have strictly positive crossing measure. 

A knot diagram is the
{\em unknot} (or is {\em unknotted}) if there is a knot $K$ having this 
diagram that is ambient isotopic to a knot $K'$ having a trivial knot 
diagram.

One can convert a knot diagram to any other diagram of an equivalent knot 
by a finite sequence
of combinatorial transformations called Reidemeister moves. 
Reidemeister showed that a knot diagram $\sK$ is unknotted if and only if 
there is a finite
sequence of Reidemeister moves that convert it to the trivial knot 
diagram. In this sense the
unknotting problem is a purely combinatorial problem, though with no 
obvious bound on the number of steps, see \cite{HL}, \cite{Wel:93b}.

\section{Unknottedness Criterion}\label{auc} 
Our approach to solving the unknotting problem, based on
that of Haken, relies on the following criterion for
unknottedness: A PL-knot $K$ embedded in $\RR^3$ is unknotted if and 
only if there exists a
piecewise-linear disk $D$ embedded in $\RR^3$ whose boundary $\partial D$ 
is the knot $K$ transversed
once. We call such a disk $D$ a {\em spanning disk}. We shall actually use 
a weaker unknottedness criterion, given
in Lemma~\ref{lem201} below.
It does not deal with a spanning disk of $K$, but rather with a spanning 
disk of another knot
$K'$ which is
ambient isotopic to $K$.

Given a PL-knot $K$, let $\sT$ be a finite triangulation of ${\bf S}^3$ 
containing $K$ in its
1-skeleton, where the 3-sphere ${\bf S}^3$ is the one-point 
compactification of $\RR^3$, and the
point ``at infinity'' is a vertex of the triangulation. Barycentrically 
subdivide $\sT$ twice to
obtain a triangulation $\sT''$, and let
$M_K = {\bf S}^3 -R_K$ denote the compact triangulated 3-manifold with 
boundary obtained by deleting
the open regular neighborhood $R_K$ of $K$. Here $R_K$ consists of all
open simplices whose closure intersects $K$. 
The closure $\bar{R}_K$ of $R_K$ is a {\em solid torus neighborhood of $K$},
a solid torus containing $K$ as its central
axis, and its boundary $\partial R_K = \partial M_K$ is homeomorphic to a 2-torus.
See \cite{Hem76} for details of this construction.
Each of $R_K$, $M_K$ and $\partial R_K = \partial M_K$
are triangulated by simplices in $\sT''$.
We call the manifold $M_K$ a {\em knot complement manifold}. 

We call a triangulation of $M_K= {\bf S}^3 - R_K$ constructed 
as above a {\em good triangulation} of $M_K$.
Similarly we define a good triangulation of a link complement manifold. 
For any good triangulation of $M_K$,
the homology group $H_1(\partial M_K; \ZZ) \approx \ZZ \oplus \ZZ $, 
since $\partial M_K$ is a 2-torus.
The kernel of the homology homomorphism induced by
the inclusion map of $\partial M_K$ into $ M_K$ is infinite cyclic (see
Theorem 3.1 of \cite{Bur:85} for a proof.)
We take as generator $(1,0)$ the 
homology class of a fixed closed oriented simple closed curve
which is the boundary $\partial B$ of an 
essential disk $B$ in $\overline{R}_K$ (a {\em meridian}) and as generator 
$(0,1)$ the homology class of a fixed closed oriented circle in $\partial 
M_K$ that has algebraic intersection $1$ with the meridian and algebraic 
linking number $0$ with $K$ (a {\em longitude}). 
A simple closed curve in $\partial \overline{R}_K$ whose homology class is
trivial in the 3-manifold
$\overline{R}_K$ but not in the surface $\partial \overline{R}_K$
is a meridian. A simple closed curve in
$\partial R_K$ whose homology class is trivial
in the 3-manifold $M_K$ but not 
in the surface $\partial \overline{R}_K$ is 
a longitude. The homology classes of a meridian and a
longitude are well-defined up to orientation.

A compact surface $S$ with boundary $\partial S$ in a compact
3-manifold with boundary $(M, \partial M)$
is said to be {\em properly embedded}
if it does not intersect itself and if $S \cap \partial M = \partial S$.
A surface $S$ is {\em essential} in $M$ if it is properly embedded in $M$, 
cannot be homotoped into $\partial M$ while holding $\partial S$ fixed, and its 
fundamental group $\pi_1 (S)$ injects into $\pi_1 (M)$ (this
last condition describes what 
topologists call an {\em incompressible surface}, see Hempel \cite{Hem76}).
In particular, a disk $S$ that is properly embedded in a 3-manifold $(M, 
\partial M)$ is essential when
$\partial S$ does not bound a disk in $\partial M$.
In the case of a knot complement $M_K$,
a properly embedded disk is essential if and 
only if the homology class $[ \partial S] \neq (0,0)$,
in $\partial M$, which happens 
if and only if removing $[ \partial S] $ does not increase the number of 
connected components of $\partial M$. 

\begin{lemma}
\label{lem201}
Let $K$ be a polygonal knot, and let $M_K$ be any good triangulation of 
${\bf S}^3 -K$.
\begin{itemize}
\item[(1)]
If $K$ is knotted, then there exists no essential disk in $M_K$.
\item[(2)]
If $K$ is unknotted, then there exists an essential disk in $M_K$. 
Furthermore any essential disk $S$ has boundary $\partial S$ representing
the homology 
class $[\partial S] = (0, \pm 1)$ in $H_1 (\partial M_K; \ZZ )$. \end{itemize}
\end{lemma}
\paragraph{Proof.}
\begin{itemize}
\item[(1)]
An essential disk can be used to give an ambient isotopy from $K$ to the
unknot.
%NOTE Explain more here?
\item[(2)]
For a knot projecting to give a trivial knot diagram, the existence of an 
essential disk in
$M_K$ is clear. For an equivalent knot $K$, the ambient isotopy carrying a 
knot with trivial
projection to $K$ carries the above disk to an essential disk for $K$. The 
homotopy class
$[\partial S ] \neq (0,0)$ since $S$ is an essential disk. The boundary of $S$ is 
a longitude,
representing the two possible generators $(0, \pm 1)$ of the kernel of the 
first homology of $\partial
M_K$ under the inclusion map of $\partial M_K$ into $M_K$.
~~~$\qed$
\end{itemize}

The Haken algorithm provides a spanning disk for the boundary $\partial S$. The 
homology condition
on $[ \partial S] $ certifies that $\partial S$ is an unoriented knot that is 
equivalent to $K$.

The usefulness of the extra condition on $[\partial S]$ is that it can be 
detected by homology with
$\ZZ / 2 \ZZ$-coefficients.

\begin{lemma} \label{lem202}
If $K$ is a PL-knot and a good triangulation $M_K = {\bf S}^3 - R_K$ 
contains a properly embedded
PL-surface S that is homeomorphic to a disk whose boundary $\partial S$ has 
homology class in $H_1( \partial M_K ; \ZZ / 2 \ZZ ) \approx \ZZ / 2 \ZZ \oplus \ZZ / 2 \ZZ$ that is not trivial,
then $K$ is
unknotted. Conversely, for any unknotted $K$ the manifold $M_K$ contains 
such a surface.
\end{lemma}
\paragraph{Proof.}
This follows from Lemma~\ref{lem201}, since $\ZZ / 2\ZZ$-homology is
obtained by reducing modulo~2. $~~~\qed$

For later use, we recall the simple fact that triangulated
surfaces that are topological disks can be recognized by their Euler 
characteristic $\chi (S) = v-e+f$, where $v,e,f$ count vertices,
edges and faces in the triangulation. 

\begin{lemma}
\label{lem203}
If $S$ is a connected triangulated surface in $\RR^3$ whose Euler 
characteristic $\chi (S) = 1$
and $\partial S \neq \emptyset$, then $S$ is homeomorphic to a disk. 
\end{lemma}
\paragraph{Proof.}
All connected surfaces have $\chi (S) \leq 2$. The only surface with $\chi 
(S) = 2$ is the sphere, which has $\partial S = \emptyset$.
The only surfaces with $\chi (S) = 1$ are the disk and the projective plane. 
The projective plane has $\partial S = \emptyset$. $~~~\qed$

\section{Normal Surfaces}\label{ns}
We work in the piecewise-linear category and let $(M, \partial M)$ be a compact 
triangulated 3-manifold with boundary.
Haken's algorithm solves the unknottedness problem by
finding a finite list of properly embedded PL-surfaces in a good 
triangulation of $M_K = {\bf S}^3 - R_K$, such that,
if $K$ is the unknot, this list will include an essential 
disk satisfying the conditions of Lemma~\ref{lem202}.
More generally, the list contains a minimal genus embedded surface bounding 
the knot. 
We will use a variation of Haken's theory based on triangulations.

A {\em normal surface} of $M$, with respect to a given triangulation, is 
a surface $S \subseteq M$ such that
\begin{enumerate}
\renewcommand{\labelenumi}{(\arabic{enumi})}
\item
$S$ is properly embedded in $M$.
\item
The intersection of $S$ with any simplex in the triangulation is transverse.
The intersection of $S$ with any tetrahedron is a 
finite disjoint union of
disks, called elementary disks. Each elementary disk is a curvilinear 
triangle or quadrilateral, whose vertices are contained on different edges 
of the
tetrahedron. A disk type is an equivalence class of elementary disks, 
where two are equivalent
if they can be deformed to one another by an isotopy preserving each 
tetrahedron.
\end{enumerate}
We are taking all disks and simplices in the above to be closed.
We allow a normal surface to have more than one connected 
component, and to be non-orientable\footnote{Some
authors require a normal surface to be connected. They call the concept
we use a {\em system of normal surfaces.}}.
Individual connected 
components may be orientable or
nonorientable surfaces.

A normal surface has associated to it combinatorial data which specifies 
the number of each
disk type that appear in the intersection of $S$ with each tetrahedron in 
the triangulation of $M$. For a given tetrahedron, each elementary
disk separates the 4 vertices into two nonempty sets; there are seven 
possibilities, consisting of
4~types of triangles which separate one vertex from the other three, and 
3~types of quadrilaterals which separate two vertices from the other 
two. If there are $t$ tetrahedra in
the triangulation of $M$ then there are $7t$ pieces of combinatorial data, 
which specify the
number of each disk type in the $t$ tetrahedra. We represent this 
combinatorial data as a
nonnegative vector
\begin{equation}
\label{eq301}
\bv = \bv (S) \in \ZZ^{7t}~,
\end{equation}
by choosing a fixed ordering of tetrahedra and of the disk types. We call 
$\bv (S)$ the {\em
normal coordinates} of $S$.

When does a vector $v \in \ZZ^{7t}$ give the normal coordinates for some normal 
surface $S$? We call
such vectors {\em admissible}. There are some obvious constraints on 
such integer admissible
vectors $\bv \in \ZZ^{7t}$.
\begin{enumerate}
\renewcommand{\labelenumi}{(\roman{enumi})}
\item
{\em Nonnegativity conditions.}
$\bv = (v_1 \dd v_{7 t} )$
has each $v_i \geq 0$.
\item
{\em Matching conditions.}
Suppose two tetrahedra $T_1 , T_2$
in the triangulation have a common face $F$. Each disk type in $T_1$ and 
$T_2$ produces either
zero or one edge in $F$ which intersects a given two of the three sides of 
$F$. The number of
edges induced between each side-pair of $F$ coming from disk types in 
$T_1$ must equal that
coming from $T_2$.
\item
{\em Quadrilateral conditions.}
In each tetrahedron in the triangulation at most one type of quadrilateral 
can occur.
\end{enumerate}
The quadrilateral conditions (iii) hold because any two quadrilaterals of 
different types placed
in a tetrahedron must intersect, which violates the embeddedness 
property of normal surfaces.

\begin{figure}[hbtp]
\centering
\includegraphics[width=.6\textwidth]{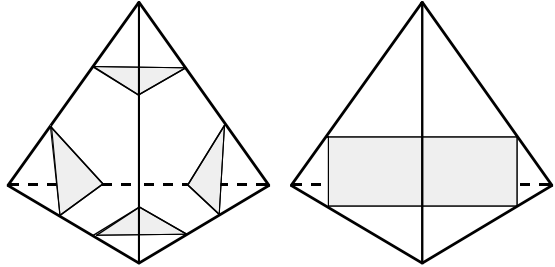}
\caption{
\label{normal}
Elementary disks in a normal surface.}
\end{figure}

We note at this point that each admissible vector $\bv = \bv (S)$ 
determines a normal surface $S$ which is unique up to a normal isotopy, 
an isotopy leaving the
triangulation of $M$ invariant. We denote this normal surface by $S( \bv )$.
Uniqueness holds, up to normal isotopy, because there is only one 
combinatorial way to
disjointly pack into a tetrahedron a given number $(n_1 , n_2 \dd n_7 )$ of 
disk types that satisfy the quadrilateral condition. The triangles
of types $1,2,3$ and 4 are stacked in parallel close to the vertex that they 
separate from
the other three vertices, all the quadrilaterals of the one type that occur 
$(5,6 ~or~ 7)$ are
arranged parallel to each other in the center of the tetrahedron.
See Figure~\ref{normal}.

The necessary conditions given above for a vector $\bv \in \ZZ^{7t}$ to be 
admissible are also sufficient.

\begin{theorem}
\label{thm301}
(Haken's Hauptsatz)
Let M be a triangulated compact 3-manifold with boundary, which consists 
of t tetrahedra. Any
integer vector $\bv \in \ZZ^{7t}$ that satisfies the nonnegativity conditions, 
matching conditions
and the quadrilateral conditions gives the normal coordinates
$\bv (S)$ of some normal surface $S$ in $M$, which is unique up to ambient 
isotopy. \end{theorem}
\paragraph{Proof.}
This is Hauptsatz 2 of Haken \cite{Hak:61}. It follows easily by noting that
the matching conditions ensure that there is a unique way to match
up the boundaries of the
elementary disks specified in each tetrahedron to form
a normal surface. $~~~\qed$ 

This result characterizes the set $\sW_M$ of all admissible vectors of 
normal surfaces as a
certain set of integer points in a rational polyhedral cone in $\RR^{7t}$. 
We define the {\em Haken normal cone}
$\sC_M$ to be the polyhedral cone in $\RR^{7t}$ cut out by the 
nonnegativity conditions and matching conditions. The points in $\sW_M$ 
are then just the integer points in the Haken normal cone $\sC_M$ that 
satisfy
the quadrilateral conditions.

Haken observed that the size of the integer vector $\bv (S)$ measures the 
{\em complexity} of the
normal surface. Define the {\em weight} $wt (S)$ of the normal surface 
$S$ to be the number
of times it intersects the 1-skeleton of $M$. 

\begin{lemma}
\label{lem301}
Suppose that $v_1 , v_2 , v_3 \in \ZZ^{7t}$ lie in Haken's normal cone 
$\sC_M$ and that
\begin{equation}
\label{eq302}
\bv_1 + \bv_2 = \bv_3~.
\end{equation}
If $\bv_3 = \bv(S_3 )$ is admissible, then so are $\bv_1$ and $\bv_2$, so 
that $\bv_1 = \bv (S_1 )$ and $\bv_2 = \bv (S_2 )$. The Euler 
characteristics of these surfaces satisfy
\begin{equation}
\label{eq303}
\chi (S) = \chi (S_1 ) + \chi (S_2 )~,
\end{equation}
and their weights satisfy
\begin{equation}
\label{eq304}
wt(S) = wt(S_1 ) + wt(S_2 ) ~.
\end{equation}
\end{lemma}
\paragraph{Proof.}
The first fact follows from Theorem~\ref{thm301}. For next two, Haken 
\cite{Hak:61} showed that $S$ is constructed from $S_1$ and $S_2$ by a 
``cutting and pasting'' operation called regular exchange,
which yields (\ref{eq303}) and (\ref{eq304}). $~~~\qed$

The ``simplest'' normal surfaces are thus those surfaces $S$ such that 
$\bv(S) \neq \bv(S_1 ) + \bv(S_2 )$
for any other (nonempty) normal surfaces $S_1$ and $S_2$. Haken calls 
these {\em fundamental surfaces},
and the corresponding vectors $\bv (S)$ {\em 
fundamental solutions}.
He proves that there are only finitely many fundamental surfaces; see 
Section~\ref{shb}. 

\begin{lemma}
\label{lem302}
If $S$ is a fundamental surface for $(M,{\partial M})$, then $S$ is connected. 
\end{lemma}
\paragraph{Proof.}
If a normal surface $S$ is disconnected, then $S = S_1 \cup S_2$ where 
$S_1$ and $S_2$ are
disjoint nonempty normal surfaces. Now $\bv (S) = \bv(S_1 ) + \bv (S_2 )$, 
contradicting $S$
being a fundamental surface.
$~~~\qed$

The basis of Haken's unknotting algorithm and knot genus algorithm is the 
following result. 

\begin{theorem}[Haken]
\label{thm302} 
Let $(M,{\partial M})$ be
a triangulated compact piecewise-linear 3-manifold $M$ with nonempty 
boundary,
which contains no embedded projective plane. If $M$ contains an essential 
surface with non-empty boundary,
then it contains such an essential surface of minimal genus, and this
surface is a fundamental normal surface.
\end{theorem}
\paragraph{Proof.}
See Chapter~5 of Haken \cite{Hak:61} or Section~5 of Schubert 
\cite{Sch:61} or \cite{JacRub89},\cite{JO:84},\cite{JT:95}.~~~~$\qed$ 

This yields:

\begin{corollary}\label{cor502}
Let $M_K$ be a knot complement manifold. There exists a properly 
embedded orientable surface $S$ in $M_K$ whose boundary is a longitude 
in $\partial M_K$ which is of minimal genus among all such surfaces, and 
which is a fundamental normal surface.
\end{corollary}
\paragraph{Proof.}
A knot complement manifold $M_K$ contains no embedded projective 
plane, since $M_K$ can be embedded in ${\bf S}^3$. Furthermore a 
minimal genus surface with boundary a longitude in $\partial M_K$ is an 
essential surface in $M_K$ \cite{Hem76}. Theorem~\ref{thm302} now 
gives an essential surface of minimal genus which is a 
fundamental normal surface.~~~$\qed$ 

In the special case where $K$ is unknotted, i.e. $K$ is of genus zero, Jaco 
and Tollefson \cite{JT:95} improved Haken's criterion given in 
Theorem~\ref{thm302}. We call a normal surface $S$ in $M$ 
a {\em vertex surface} if $\bv (S)$ lies on an
extremal ray (1-dimensional face) of the Haken normal cone $\sC_M$. In 
this case we call $\bv (S)$ a {\em vertex solution} of $\sC_M$.
The notion of a vertex surface 
was introduced by Jaco and Oertel \cite{JO:84}, who imposed the 
addtional requirement that it be a connected orientable surface that is 
either a fundamental surface or twice a fundamental surface.
The last case can occur when a surface $S$ is the boundary of 
a regular neighborhood of another
surface which is one-sided. However
$\chi(S)$ is even in this case, so that this case cannot
occure when the surface is a disk. 

\begin{theorem} [Vertex Surface Theorem]
\label{thm303}
If a triangulated compact 3-manifold $M$ with nonempty boundary $\partial M$ 
contains an essential
disk, then it contains such a disk which is a vertex surface.
\end{theorem}
\paragraph{Proof.}
This is an immediate
 consequence of  Corollary 6.4 of Jaco and Tollefson \cite{JT:95}. 
$~~~\qed$

The theory of normal surfaces can also be applied to determine the 
splittability of links.

\begin{theorem} [Splittability Criterion]
\label{thm504}
A link $L$ is splittable if and only if any link complement manifold $M_L$ 
contains a fundamental normal 2-sphere $S$ that separates two 
components of $\partial M_L$.
\end{theorem}
\paragraph{Proof.}
This is established in Section 4 of Schubert \cite{Sch:61}, who attributes 
the result to Haken \cite{Hak:61}. See also \cite{JacRub89} \cite{JT:95}~~~$\qed$

This result was strengthened by Jaco and Tollefson \cite{JT:95}. 
\begin{theorem}\label{thm505}
If a link complement manifold $M_L$ contains a normal 2-sphere $S$ which 
separates two components of $\partial M_L$, then it contains such a sphere 
that is a vertex surface. \end{theorem}
\paragraph{Proof.} This follows from Theorem~5.2 \cite{JT:95}.~~~$\qed$ 

To check that a sphere $S$ separates two components of $\partial M_L$ it 
suffices to exhibit a path in the 1-skeleton of $M_L$ that connects two 
components of $\partial M_L$ and whose intersection number with $S$ is odd. 

For the complexity analysis, we must give bounds for the number of 
vectors $\bv$ in
the exhausting set, and for the sizes $||\bv ||_1$ of these vectors, where 
the size is
\begin{equation}
\label{eq305} || \bv ||_1 := \sum_{i=1}^{7t} v_i ~. \end{equation}
We address these questions in the next section. Then in Section~\ref{tri} 
we deal with the problem of constructing the triangulated 3-manifold 
$M_K = {\bf S}^3 -R_K$ from a knot diagram $\sK$. Sections 8--10 give 
the complexity
analysis.

\section{Bounds for Fundamental Solutions and Hilbert Bases}\label{shb} 
We bound the number and size of fundamental solutions in the Haken 
normal cone $\sC_M$
of an arbitrary triangulated compact 3-manifold $M$ with boundary $\partial 
M$ that contains $t$
tetrahedra. The system of linear inequalities defining Haken's normal cone 
$\sC_M$ has the
form:
\begin{equation}
\label{eq401}
v_i \geq 0 \quad (1 \leq i \leq 7t)~,
\end{equation}
\begin{equation}
\label{eq402}
v_{k_1} + v_{k_2} = v_{k_3} + v_{k_4} \quad \mbox{(at most}~ 6t ~ 
\mbox{equations}) ~. \end{equation}
The matching conditions (\ref{eq402}) all involve exactly four 
variables as indicated,
because only two of the seven disk types in a tetrahedron $T$ have edges 
parallel to one specified
edge on one face of $T$. There are at most $6t$ such equations because the 
$t$ tetrahedra have
between them $12t$ such edge pairs, and each equation matches two pairs 
that occur in no other
equation. (There will be strictly less than $6t$ equations if $M$ has 
nonempty boundary.) The
cone
$\sC_M$ is a {\em pointed cone}, (i.e. it contains no line) because it is 
contained in the
positive orthant $\RR^{7t}_+$. It is not full-dimensional because it has 
equality constraints,
but
(\ref{eq402}) implies
\begin{equation}
\label{eq403}
t \leq \dim_\RR ( \sC_M ) \leq 7t~.
\end{equation}
It is a {\em rational cone} because it is cut out by rational equality and 
inequality
constraints; equivalently, each of its extreme rays contains an integral 
vector.
The concept of fundamental solutions of $\sC_M$ is closely related to 
that of a Hilbert basis for
the rational cone $\sC_M$.
Given a (homogeneous) rational cone $\sC$ in $\RR^m$, its set of 
integral points
\begin{equation}
\label{eq404}
\sC( \ZZ ) := \sC \cap \ZZ^m
\end{equation}
forms a semigroup under addition.
An (integral)
{\em Hilbert generating set}
$\sG$ for $\sC$ is any finite set of generators for $\sC( \ZZ )$. If 
$\sC$ is a pointed rational cone, then $\sC$ has a unique minimal 
generating set $\sH$, see
Schrijver\cite[Theorem~16.4]{Sch86}.
It consists exactly of all elements $\bv \in \sC ( \ZZ )$ which are {\em 
minimal} in the sense that
\begin{equation}
\label{eq405}
\bv \neq \bv_1 + \bv_2 \quad \mbox{for any nonzero}~~ \bv_1 , \bv_2 \in 
\sC( \ZZ) ~.
\end{equation}
We call this set
$\sH = \sH ( \sC )$ the
{\em minimal Hilbert basis}
of $\sC$.
By definition, any fundamental solution of $\sC_M$ is in the minimal 
Hilbert basis $\sH ( \sC_M )$.
There may, however, be elements in $\sH ( \sC_M )$ that are not 
fundamental solutions because they violate the quadrilateral conditions.

A {\em minimal vertex solution} of a pointed homogeneous rational cone 
$\sC$ in $\RR^m$ is the smallest nonzero integral point on any extreme 
ray (1-dimensional face) of the cone $\sC$. A minimal vertex solution 
is always in the minimal Hilbert basis $\sH ( \sC )$, because
the only points in $\sC$ that it
can be a convex combination
of are other points in the
ray $\RR^+ \bv = \{ \la \bv : \la \geq 0 \}$. By definition a fundamental 
vertex solution of the Haken cone $\sC_M$ is a minimal vertex solution 
of $\sC_M$.

\begin{lemma} \label{lem401}
Let $M$ be a triangulated compact 3-manifold, possibly with boundary, that 
contains $t$ tetrahedra in the triangulation.
\begin{itemize}
\item[(1)]
Any vertex minimal solution $\bv \in \ZZ^{7t}$ of the Haken normal cone 
$\sC_M $ in $\RR^{7t}$ has
\begin{equation}
\label{eq406} \max_{1 \leq i \leq 7t} (v_i ) 
\leq 2^{7t-1} ~. \end{equation}
\item[(2)]
Any minimal Hilbert basis element $\bv \in \ZZ^{7t}$ of the Haken normal 
cone $\sC_M$ has
\begin{equation}
\label{eq407}
\max_{1 \leq i \leq 7t} (v_i ) \leq t \cdot 2^{7t+2} ~. \end{equation}
\end{itemize}
\end{lemma}
\paragraph{Proof.}
\begin{itemize}
\item[(1)]
Choose a maximal linearly independent subset of matching conditions
(\ref{eq402}). There will
be $7t-d$ of them, where
$d = \dim_\RR ( \sC_M )$.
Any vertex ray is determined by adjoining to these equations $d-1$ other 
binding inequality constraints
\begin{equation}
\label{eq408} \{ \bv_{i_k} = 0 ~; \quad 1 
\leq k \leq d-1 \}~, \end{equation} with the proviso that the resulting system have 
rank $7t-1$. These conditions yield a
$(7t-1) \times 7t$ integer matrix $\bM$ of rank $7t-1$, and the vertex 
ray elements $\bz = (z_1 \dd z_{7t} )$
satisfy
\begin{equation}
\label{eq409}
\bM \bz = {\bf 0} ~.
\end{equation}
In order to get a feasible vertex ray in $\sC_M$ all nonzero coordinates 
$z_i$ must have the same sign.

At least one of the unit coordinate vectors $\be_1 , \be_2 \dd \be_{7t}$ 
must be linearly independent
of the row space of $\bM$. Adjoin it as a first row to $\bM$ and we 
obtain a full rank $7t \times 7t$ integer matrix $$
\tilde{\bM} =
[ \begin{array}{c}
\be_k \\
\bM
\end{array} ] ~, \quad \det ( \tilde{\bM} ) \neq 0 ~. $$ Consider the 
adjoint matrix
$ adj ( \tilde{\bM}) = \det ( \tilde{\bM}) \tilde{\bM}^{-1}$, which 
has integer entries
\begin{equation}
\label{eq410}
w_{ij} : = adj ( \tilde{\bM})_{ij} = (-1)^{i+j} \det \tilde{\bM} [j|i ]~,
\end{equation}
in which
$\tilde{\bM} [j|i ]$ is the
minor obtained by crossing out the
$j^{\rm th}$ row and $i^{\rm th}$ column of $\tilde{\bM}$. Let 
\begin{equation}
\label{eq411}
\bw = [w_{11}, w_{21} \dd w_{7t,1} ]^t
\end{equation}
be the first column of $adj ( \tilde{\bM} )$. Since $\tilde{\bM} adj ( 
\tilde{\bM} ) = \det ( \tilde{\bM} ) \bM$ this yields $$ \tilde{\bM} 
\bw = {\bf 0 }~,
$$
and $\bw \neq {\bf 0}$ because $\langle \be_k , \bw \rangle = \det ( 
\tilde{\bM} ) \neq 0$. We bound the entries of $\bw$ using Hadamard's 
inequality, which states for an $m \times m$ real matrix that $$ \det ( 
\bN )^2 \leq \prod_{i=1}^m || \bn_i ||^2 ~, $$ in which $|| \bn_i||^2$ is the 
Euclidean length of the $i^{\rm th}$ row $\bn_i$ of $\bN$.
We apply this to (\ref{eq409}), and observe that each row of the $(7t-1) 
\times (7t-1)$ matrix $\tilde{\bM} [j|i ]$ has squared Euclidean length 
at most 4, because this is true for all row vectors in the system 
given by (\ref{eq402}) and for $\be_k$.
Applied to (\ref{eq409}) this gives
$$
|w_{ij} |^2 \leq 4^{7t-1} ~.
$$
However (\ref{eq409}) shows that $\bw \in \ZZ^{7t}$, and a vertex minimal 
solution $\bv$ in the extreme ray is obtained by dividing $\bw$ by the 
greatest common divisor of its elements, hence $$
\max_{1 \leq i \leq 7t} |v_i | \leq \max_{1 \leq i \leq 7t} |w_{i1} | \leq 
2^{7t-1} ~.
$$
\item[(2)]
A {\em simplicial cone} $\sC$ in $\RR^{7t}$ is a $d$-dimensional 
pointed cone which has exactly $d$ extreme rays. Let $\bv^{(1)} \dd 
\bv^{(d)} \in \ZZ^{7t}$ be the vertex minimal solutions for the extreme 
rays. Each point in $\sC$ can be expressed as a nonnegative linear 
combination of the $\bv^{(j)}$, as $$ \bv = \sum_{j=1}^d \la_j \bv^{(j)} ~, 
\quad \mbox{each}~~ \la_j \geq 0 ~. $$
If $\bv$ is in the minimal Hilbert basis $\sH ( \sC )$, then $0 \leq \la_j 
\leq 1$, for otherwise one has $$
\bv = ( \bv - \bv^{(j)} ) + \bv^{(j)} ,
$$
and both $\bv - \bv^{(j)}$ and $\bv^{(j)}$ are nonzero integer vectors in 
$\sC$, which is a contradiction.
Thus,
any minimal Hilbert basis element $\bv = (v_1 , \ldots , v_{7t} )$ of a 
simplicial cone satisfies
\begin{equation}
\label{eq412}
|v_i | \leq \sum_{j=1}^d |v_i^{(j)} | \quad \mbox{for}~~ 1 \leq i \leq 7t~. 
\end{equation}

The cone $\sC_M$ may not be simplicial, but we can partition it into a 
set of simplicial cones
$\{ \sC_M^{(k)} \}$ each of whose extreme rays are extreme rays of 
$\sC_M$ itself.
We have
$$
\sH ( \sC_M ) \subseteq \bigcup_{k} \sH ( \sC_M^{(k)} ) ~. $$ Thus all 
Hilbert basis elements of $\sH ( \sC_M )$ satisfy by the bound 
in (\ref{eq412}) for Hilbert basis elements of $\sH ( \sC_M^{(k)} )$. Using 
(\ref{eq406}) to bound $|v_i^{(j)} |$ we obtain $$ |v_i | \leq d 2^{7t-1} \leq 
7t~2^{7t-1} \leq t ~2^{7t+2} ~, $$ as required.
$~~~\qed$
\end{itemize}	

\begin{lemma}
\label{lem402}
\begin{itemize}
\item[(1)]
The Haken normal cone $\sC_M$ has at most $2^{7t}$ vertex 
fundamental solutions.
\item[(2)]
The Haken normal cone $\sC_M$ has at most $t^{7t} ~ 2^{49t^2 + 14t}$ 
elements in its minimal Hilbert basis.
\end{itemize}
\end{lemma}

\paragraph{Proof.}
\begin{itemize}
\item[(1)]
There are $\binom{7t}{d-1}$
possible choices for the systems given in (\ref{eq408}),
and ${\binom{7t}{d-1}} \le 2^{7t}$.
\item[(2)]
We wastefully count every integer vector $\bv$ in the box $0 \leq v_i \leq t 
2^{7t+2}$ for $1 \leq i \leq 7t$ that is allowed by Lemma~\ref{lem401}. $~~~\qed$
\end{itemize}
\paragraph{Remark.}
The bounds in Lemma~\ref{lem402} are very crude. However the functional 
forms $2^{c_1 t}$ and $2^{c_2
t^2}$ of the bounds are best possible for general rational cones, apart 
from the values of $c_1$ and
$c_2$.

\section{Triangulations}\label{tri}
Given a link diagram $\sD$ with crossing measure $n$, we show how to 
construct a triangulated
3-manifold
$M_L \cong {\bf S}^3 -R_L$, where $R_L$ is a regular neighborhood of a 
link $L$ which has a regular
projection that is the link diagram $\sD$. The construction takes time 
polynomial in the crossing measure of
$\sD$, and the triangulations of each of $M_L$ and $\bar{R}_L$ contain 
$O(n)$ tetrahedra.

\begin{lemma}
\label{lem501}
Given a link diagram $\sD$ of crossing measure $n$, one can construct in 
time $O(n \log n )$ a
triangulated convex polyhedron P in $\RR^3$ such that: \begin{enumerate}
\renewcommand{\labelenumi}{(\roman{enumi})}
\item
The triangulation has at most $420n$ tetrahedra. \item
Every vertex in the triangulation is a lattice point $(x,y,z) \in \ZZ^3$, with 
$0 \leq x \leq 30n$, $0 \leq y \leq 30n$ and $-4 \leq z \leq 4$. \item There is a link L 
embedded in the
1-skeleton of the triangulation which lies entirely in the interior of $P$, 
and whose orthogonal
projection on the $(x,y)$-plane is regular and is a link diagram isomorphic 
to $\sD$.
\end{enumerate}
\end{lemma}
\paragraph{Proof.}
We first add extra vertices to $\sD$.
To each non-loop edge we add a new vertex of degree~2, which splits it 
into two edges. To
each non-isolated loop we add two new vertices of degree~2, which split 
it into three edges. To each isolated loops we add three new vertices of 
degree~2, making it an isolated
triangle. The resulting labelled graph $\sG$ is still planar, has no loops or 
multiple edges,
and has at most $5n$ vertices. (The worst case consists of several 
disjoint single
crossing projections.) Let $m$ denote the number of vertices of $\sG$, and 
call the $m-n$
vertices added {\em special vertices}.

Using the Hopcroft-Tarjan planarity testing algorithm \cite{HopTar74} we 
can construct a planar embedding of $\sG$ in time $O(n \log n )$. 
From this
data we determine the planar faces of this embedding, and add extra edges 
to triangulate each face, thus obtaining a triangulated planar graph $\sG'$,
in time $O(n)$. The graph $\sG'$ has $m$ vertices and $2m-5$ bounded 
triangular faces, and the unbounded face is also a triangle.

It was shown by de Frajsseix, Pach and Pollack \cite{deFPacPol90} that 
there exists a planar embedding of $\sG'$ whose vertices
$v( \sG' )$ lie in the plane $z = -1$, in the grid
\begin{equation}
\label{eq502}
\{(x,y,-1): 0 \leq x,y \leq 10n-1 ~; \quad x,y \in \ZZ \}~, \end{equation}
in which all edges of $\sG'$ are straight line segments. 
They also show in \cite[Section~4]{deFPacPol90} that one can explicitly 
find such an embedding in time $O( n \log n )$.

Next we make an identical copy $\sG''$ of $\sG'$ on the plane $z=1$ with 
$\sG'' = \sG + (0,0,2)$, and
vertex set $V( \sG'') = V( \sG ) + (0,0,2)$. We now consider the polyhedron 
$P'$ which is the convex
hull of $V( \sG')$ and $V( \sG'')$. It is a triangular prism, because the 
outside face of $\sG'$ is
a triangle. We add $m_1$ vertical edges connecting each vertex $\bv \in 
V( \sG' )$ to its copy
$\bv + (0,0,1) \in V( \sG'')$.
Let $\sE$ denote these edges, together with all the edges in $\sG'$ and 
$\sG''$. Using these edges, the polyhedron $P'$ decomposes into $2m-5$ 
triangular prisms $\{Q_j : 1 \leq
j \leq 2m-5\}$, with top and bottom faces of each being congruent 
triangular faces of $\sG'$ and
$\sG''$.

We next triangulate $P'$ by dissecting each of the $2m-5$ triangular 
prisms into 14
tetrahedra, as follows. We subdivide each vertical rectangular face into 
four triangles using its diagonals. Then we cone each rectangular face to 
the centroid of $P'$, and note that the
centroid lies in the plane $z = 0$. We add 4 new vertices to each prism, 
one on each face and
one in the center. The point of this subdivision is that the triangulations 
of adjacent prisms
are compatible. Note also that all new vertices added lie in the lattice 
$\displaystyle \frac{1}{6} \ZZ^2$.

Let $\sE'$ denote all the edges in the union of these triangulated prisms. 
We identify the link diagram $\sD$ with the graph $\sG$ embedded in 
$\sG'$. We next observe that
there is a polygonal link $L$ imbedded in the edge set $\sE'$ whose 
projection is the link
diagram $\sD$. We insist that any edge in $\sD$ that runs to an 
undercrossing have its
undercrossing vertex lie in the plane $z = -1$, while any edge that runs to 
an overcrossing has
its overcrossing vertex lie in the plane $z = 1$. Each such edge travels 
from one of the
$n$ original vertices
of $\sD$ to a special vertex. Edges that do not meet vertices labelled 
overcrossing or
undercrossing are assigned to the $z = -1$ plane. The edge corresponding 
to an edge running from
the $z = -1$ to the $z = 1$ plane is contained in one of the diagonals added 
to a
prism. The resulting embedding $L$ in $P'$ has a regular orthogonal 
projection onto the
$(x,y)$-plane, since no edges are vertical and only transverse double 
points occur as singularities
in the projection.

The knot edges lie in the boundary of $P'$. To correct this we take two 
additional copies of
$P'$ and glue one to its top, along the plane $z =1$, and one to its bottom, 
along the plane $z
=-1$. The knot now lies in the interior of the resulting polyhedron $P''$. 
All the vertices added
in this construction lie in $\displaystyle \frac{1}{6} \ZZ^3$, so a
rescaling by a factor of 6 puts all vertices in $\ZZ^3$.
The total number of tetrahedra used in
triangulating $P''$ is $84m$, which is at most $420n$. 

Now $P''$ satisfies properties (i), (ii) and (iii).$~~~\qed$ 

We now construct a triangulated 3-manifold $M_L$ embedded in ${\bf S}^3$
from the triangulation of $P$ in Lemma~\ref{lem301}.
This construction only
uses the {\em combinatorial triangulation} of $P$, i.e. the adjacency 
relations of tetrahedra and of the set of edges
of $P$ that specify the link $L$, and does not use its given embedding in 
$\RR ^3$.

\begin{lemma}
\label{lem502}
Given a link diagram $\sD$ of crossing measure n, one can construct in 
time $O( n \log n)$ a
combinatorial triangulation of ${\bf S}^3$ using at most $253440(n+1)$ 
tetrahedra, which contains a good
triangulation of $M_L \cong {\bf S}^3 -R_L$, with $\bar{R}_L$ a regular 
neighborhood of a combinatorial
link $L$ with link diagram $\sD$, and $\partial R_L = \partial M_K$.
Furthermore one 
can construct in time $O(n^2 \log n)$ in the triangulation marked sets of 
edges in $\partial M_L$ for a meridian on each 2-torus 
component of $\partial M_L$, and a set of marked paths of $O(n)$ edges in $M_L 
\setminus \partial M_L$ that connect pairs of components of $\partial M_L$, which 
between them connect all components of $\partial M_L$.
\end{lemma}

\paragraph{Proof.}
We construct a combinatorial triangulation $\sT$ of ${\bf S}^3$ by coning 
all the triangular outside faces of the triangulation of $P$ to the ``point at 
infinity'' added in
constructing ${\bf S}^3$ as the one point compactification of $\RR^3$. 
This adds $4m+2$ new tetrahedra, where $m \leq 5n$, for a total of at most 
$440n+2$ tetrahedra. Barycentrically subdivide twice to
obtain a triangulation $\sT''$, and take $M_L$ to be ${\bf S}^3$ minus the 
interior of a regular
neighborhood of $L$. Since barycentric subdivision splits a tetrahedron 
into 24~tetrahedra, we use at
most $24^2 (440n + 2)$ tetrahedra. It is easy to determine a meridian on 
the 2-torus
in $\partial M_L$ corresponding to each component of $L$, in time $O(n)$. To 
construct a longitude we use the projection of the polyhedron $P$ in 
$\RR^3$ on the z-axis in Lemma~\ref{lem501} to construct a path with 
algebraic linking number zero with the core of the 2-torus. This can be 
done by suitably making a ``twist'' at each vertex. 

We can find a shortest path between each pair of components of $\partial M_K$ 
in time $O( n \log n )$ each. We retain a minimal set of such paths which 
contain no edge inside any component of $\partial M_K$,
whose union connects all components of $\partial M_K$. If done 
carefully, this can be done in time $O(n^2 \log n)$ and result in a list of at 
most $O(n)$ paths involving $O(n^2)$ edges.~~~$\qed$ 

\section{Certifying Unknottedness}\label{cuua} 
To show that the {\sc unknottedness problem} is in {\bf NP}, we must 
construct for any $n$-crossing knot diagram $\sD$ a polynomial length 
certificate that proves it is unknotted. The certificate takes the following 
form. 

\paragraph{Unknottedness Certificate}
\begin{enumerate}
\item
Given a link diagram $\sD$ with $n$ crossings certify that $\sD$ is a knot 
diagram.
\item
Construct a piecewise-linear knot $K$ in $\RR^3$ which has regular 
projection $\sD$. From it construct a
good triangulation $M_K \cong {\bf S}^3 - \bar{R}_K$ which contains $t$ 
tetrahedra, with $t = O(n)$, and
with a meridian of $\partial R_K$ marked in $\partial M_K$. \item Guess a suitable 
fundamental solution $\bv \in \ZZ^{7t}$ to the Haken normal equations for 
$M_K$.
Verify the Haken quadrilateral conditions. Let $S$ denote the associated 
normal surface, so $\bv= \bv (S)$.
\item Certify that $S$ is an essential disk.
\begin{description}
\item{(4a)}
Certify that $S$ is connected.
\item{(4b)}
Certify that $S$ is a disk.
Check that $\partial S \neq \emptyset$ and that it
has Euler characteristic $\chi (S) = 1$.
\item{(4c)}
Certify that $S$ is essential by checking that its homology class $[ \partial S] 
\neq (0,0)$ in $H_1 ( \partial M_K ; \ZZ / 2 \ZZ )$.
\end{description}
\end{enumerate}

The approach of Haken to deciding unknottedness is based on an exhaustive 
search for a certificate of the above kind. 

\section*{Haken-type Unknottedness Algorithm} \begin{description} 
\item{\em Input:}
Link diagram $\sD$ of crossing number $n$. \item{\em Question:} Is $\sD$ 
a knot diagram of the unknot?
\end{description}

\begin{enumerate}
\item
Test if $\sD$ is a knot diagram.
If so, denote it $\sK$.
\item
Construct a finite combinatorial triangulation $\sT$ of ${\bf S}^3$ which 
contains a polygonal knot $K$ in
its one-skeleton which has diagram $\sK$, and also contains a good 
compact triangulated 3-manifold
$M_K \cong {\bf S}^3 -R_K$,
where $\bar{R}_K$ is a regular neighborhood of $K$. Let $t$ denote the 
number of tetrahedra in
$M_K$. \item Construct an exhausting list $\sL$ of vectors $\sL 
\subseteq \ZZ^{7t}$ that
contains $\bv (S)$ for some embedded compressing disk in $M_K$, if one 
exists. \item
For each $\bv \in \sL$, test if $\bv$ is admissible. If so let $S$ denote 
the normal surface with $\bv = \bv (S)$. Test if $S$ is an essential disk 
for $M_K$. 

\begin{description}
\item{(4a)}
Is $S$ connected?
\item{(4b)}
Is $S$ a disk?
Check that the Euler characteristic $\chi (S) = 1$ and $\partial S \neq 
\emptyset$. \item{(4c)}
Is $S$ essential?
Check that the homology class $[ \partial S ] \in H_1( \partial M; \ZZ/ 2 \ZZ )$ is 
nontrivial by computing that
its intersection number $(\bmod ~2)$ with a meridian of the 2-torus 
$T_K$ is $ 1 \ ( \bmod ~2)$.
\end{description}
If all tests are passed for $S$,
answer ``yes,'' halt and output $\bv (S)$. 

\item
If the complete list $\sL$ is examined and no essential disk is found, 
answer, ``no'' and halt.
\end{enumerate}

For the
{\em Haken unknottedness algorithm} we take for the Exhausting List 
$\sL_H : = \sL_H ( M_K )$ the set
\begin{equation}
\label{eq601}
\sL_H : = \{ \bv = (v_1 \dd v_{7t} )
\in \ZZ^{7t} : 0 \leq v_i \leq t2^{7t+2} \}~. \end{equation} For the {\em Jaco-Tollefson 
unknottedness algorithm} we take for the exhausting list
$\sL_V = \sL_{\bv} (M_K)$ the set
\begin{equation}
\label{eq602}
\sL_{V} : = \{ \bv = (v_1 \dd v_{7t} ) \in \ZZ^{7t} : \bv ~~\mbox{a 
minimal vertex solution in}~~
\sC_M \}~.
\end{equation}

\begin{theorem}
\label{Nth81}
\begin{itemize}
\item[(1)]
There is a constant $c$ such that the Haken unknottedness algorithm 
decides for any $n$-crossing link diagram whether it is a knot diagram 
that represents the trivial knot using at most $O( \exp (cn^2))$ time and 
$O(n^2 \log n )$ space, on a Turing machine.
\item[(2)]
There is a constant $c'$ such that the Jaco-Tollefson unknottedness 
algorithm decides the same question in at most time $O( \exp (c'n ))$ time 
and $O(n^2 \log n )$ space, on a Turing machine. \end{itemize}
\end{theorem}

\paragraph{Proof.}
Assuming that the individual steps of the algorithms 
are proved correct, the two algorithms are correct by
Lemmas~\ref{lem201} and \ref{lem202}, respectively. \\ 

(1)~For the Haken unknottedness algorithm, the list $\sL_H$ is 
exhausting
by Theorem~\ref{thm302} with Lemma~\ref{lem201}. 

Step (1) is done in polynomial time $O(n \log n )$ by tracing a path through 
the link diagram $\sD$ and
checking that it visits every edge of $\sD$. 

Step (2) is done in polynomial time $O(n \log n )$ by 
Lemmas~\ref{lem501} and \ref{lem502}. The
triangulated manifold $M_K$ has at most
$253440(n+1)$ tetrahedra.

The list $\sL_H$ is a large cube containing $t^{7t} 2^{49t^2 +14t}$ 
integer vectors, using integers containing at most $7t+ \log_2 t +2$ 
binary digits.
We now test each $\bv \in \sL_H$, one at a time, to determine whether 
it is a fundamental solution that is an essential disk.

To begin step (4), given
$\bv \in \sL_H$,
we first test whether $\bv$ corresponds to a normal surface by checking 
whether the quadrilateral
conditions hold: for each tetrahedron in $M_K$ at least two of the three 
quadrilateral variables
vanish. This takes time $O(t)$ and space $O(t^2 )$. Now we have $\bv = 
\bv(S)$ for some normal
surface $S$. We can treat $S$ as a triangulated surface by splitting 
quadrilaterals into
four triangles by adding two diagonals to each quadrilateral.

If $S$ is connected, then we can check whether $S$ is homeomorphic to a 
disk in polynomial time
$O(n^2 \log n )$.
We use the criterion of Lemma~\ref{lem203}. We check that $\partial S \neq 
\emptyset$ by checking that $v_i \neq 0$ for some 
triangle or quadrilateral is that has an edge in $\partial M$.
We can compute the 
Euler characteristic $\chi (S)$ in polynomial time directly from the data 
in $\bv$, since it encodes the necessary data on the vertices, edges and 
faces of the triangles and quadrilaterals in $S$.
According to Jaco and Tollefson \cite[p.~404]{JT:95} we have
\begin{equation}
\label{eq603}
\chi (S) = \df{1}{2} S_3 - \sigma (S) + wt (S) \end{equation} in which $S_3$ counts 
the total number of triangles in $S$, 
$\sigma (S) = \sum_{i=1}^{7t} v_i$, and wt(S) is defined below. For each 
edge $j$ in the triangulation of $M_K$, let $t_j$ count the number of 
tetrahedra
that contain edge $j$, and set $\eps_{ji} = 1$ if edge $e_j$ meets a disk 
of type $i$.
Then
\begin{equation}
\label{eq604}
wt(S) = \sum_{i=1}^{7t} \sum_j \df{\eps_{ji} v_i}{t_j} ~. \end{equation} 
The weight counts the intersection of $S$ with the 1-skeleton
of the triangulation. Thus $\chi (S)$
is computed by (\ref{eq603}) using $O(n^2 \log n )$ bit
operations, and we can check whether $\chi (S) = 1$.

If $S$ is known to be a disk, then we can check whether $\partial S$ is an 
essential curve in $\partial M$ in polynomial
time. Since $S$ is a disk the boundary $\partial S$ must have homology class 
$(0,0)$ or $(0, \pm 1)$
in $H_1 ( \partial M ;\ZZ)$,
by Lemma~\ref{lem201}.
We distinguish these possibilities by determining if the class $[ \partial S] \in 
H_1 ( \partial M_1 ; \ZZ / 2 \ZZ )$ is $(0,0)$ or $(0,1)$.
To do this we compute the intersection number $( \bmod ~2)$ of $\partial S$ 
with a fixed meridian in the
1-skeleton of the 2-torus $\partial M_K$. Such a meridian $\ga$ can be 
determined in time $O( n \log n )$ in the
process of constructing $M_K$ in step (2). This intersection number can be 
computed in time $O(n \log n )$ as
\begin{equation}
\label{eq605}
n_S : = \df{1}{2} ( \sum v_i ) ~~( \bmod ~2 )~, \end{equation} in which the 
sum is taken over those disk types $i$ that correspond to triangles
and quadrilaterals with an edge on $\partial M_K$ and a vertex on
that edge that is contained in some edge of $\ga $.
We count disk types $i$ with
multiplicity equal to the number of vertices of this type contained in the
disk type. Then $S$ is essential only
if $n_s \equiv 1 ( \bmod ~2)$, using Lemma~\ref{lem202}. 

Thus to test a given $\bv \in \sL_H$
requires $O(t^{7t} 2^{49t^2 + 14t} )$ operations, and testing $S$ for being a 
fundamental solution is
the only step that requires exponential time. 

The total running time of the Haken algorithm is dominated by the 
test for being a fundamental
solution and the size of the
exhausting list $\sL_H$. It is easy to sequentially test elements of 
$\sL_H$ in lexicographic
order, so that the space complexity does not materially increase over that 
for testing a single
vector $\bv$. Since $t \leq O(n)$ we obtain upper bounds of the form 
$O(2^{cn^2} )$ for time
complexity and $O(n^2 \log n )$ for the space complexity. 

(2).
For the Jaco-Tollefson algorithm, the list $\sL_{V}$ is exhausting by
the Vertex Surface Theorem~\ref{thm303} and Lemma~\ref{lem401}. 
There are at most $2^{7t}$
elements in $\sL_{V}$
by Lemma~\ref{lem402}. Furthermore these elements are easy to 
enumerate sequentially in time
$O( \log n )$ each:
First precompute a maximal linearly independent set of $7t-d$ matching 
condition constraints,
where $d = \dim_\RR ( \sC_M )$.
Then test all $(d-1)$-duples
$\{i_1 \dd i_{d-1} \} \subseteq \{1,2 \dd 7t \}$ to see if adjoining the 
constraints
$$
x_{i_j} = 0 ~, \quad 1 \leq j \leq d
$$
gives an extremal ray of the Haken cone $\sC_M$. If so, determine a 
minimal point $\bv \in
\ZZ^{7t}$ on this ray and run the rest of the algorithm. Given $\{i_1 \dd 
i_j \}$ one
can compute $\bv$ in polynomial time
$O(n^2 \log n )$
by first computing $\bw$ in
(\ref{eq410}) and dividing it
by g.c.d.
$(w_1 , w_2 \dd w_{7t} )$ to get $\bv$.

The testing time, for each $\bv \in \sL_{V}$, for whether $\bv = \bv(S)$ 
gives an essential disk
$S$ goes down to $O(t2^{7t} )$ time and to $O(n^2 \log n )$ space, because 
the
computationally intensive
test whether $S$ is connected is eliminated: $S$ is connected because 
minimal vertex solutions are
always Hilbert basis elements, so $S$ is a fundamental surface and 
connected by Lemma~\ref{lem301}. 

Since $t = O(n)$ we obtain the time bound $O(2^{cn} )$ and space bound at 
most $O(n^2 \log n )$.

We check whether $S$ is fundamental by exhaustive search. For each $\bw 
\in \ZZ^{7t}$ with
$0 \leq w_i \leq v_i$ for $1 \leq i \leq 7t$ we set $\bw' = \bv - \bw$ and 
check if
$\bw , \bw' \in \sC_M$.
If this happens for some $\bw$ then
$\bv = \bw + \bw'$ is not fundamental.
If this happens for no $\bw$ then $\bv$ is in the Hilbert basis, and is a 
fundamental surface by
Lemma~\ref{lem301}.
There are at most $O(t^{7t} 2^{49t^2 + 14t} )$ values of $\bw$ to check. We 
can check them sequentially, using space $O(t^2 )$. If $S$ is a fundamental 
surface, then it is
guaranteed to be connected by Lemma~\ref{lem302}.~~~$\qed$ 

We can now prove that the Unknottedness Certificate above can be 
checked in polynomial time.

\paragraph{Proof of Theorem~\ref{knotnp}}
If all the steps in the certificate are checked then $\sD$ is unknotted by 
Lemma~\ref{lem202}. It remains to show
each step can be done in length polynomial in $n$.
\begin{enumerate}
\item
This is trivially polynomial time in $n$, by tracing edges around the graph.
\item
This is polynomial time by Lemma~\ref{lem502}, upon observing that $t = 
O(n)$.
\item
This step is nondeterministic. According to Lemma~\ref{lem401},
any fundamental solution is
$$
wt( \bv ) \leq \sum_{i=0}^{7t} v_i \leq 7t^2 2^{7t+2} ~,
$$
hence $\bv$ can 
be guessed in nondeterministic time linear in $t$, which is $O(n)$.
At this point we do not
have a proof that the solution is fundamental. Verifying the quadrilateral 
conditions takes time
$O(n^2 )$, since we have $O(n)$ tetrahedra and have $O(n)$-bit integers
in $\bv$ to examine.
\end{enumerate}
\begin{description}
\item{(4a)}
According to the Vertex Surface Theorem (\ref{thm303}), we may choose a 
vertex fundamental solution in (3)
for which $S$ is an essential disk. We may certify in polynomial time that
$\bv = \bv (S)$ is in the Hilbert basis, by guessing the correct linearly 
independent set of $7t-1$ constraints that are binding for $\bv$, and then 
verifying that $\bv \in \ZZ^{7t}$ is primitive, i.e.
$$
g.c.d. (v_1 , v_2 \dd v_{7t} ) = 1~.
$$
This takes time $O(n^3 \log n )$.
We verify the $7t-1$ constraints are linearly independent by guessing 
which one of the unit coordinate vectors $\be_1 \dd \be_{7t}$ extends to 
a basis of $\RR^{7t}$ and verifying that the
determinant of the resulting matrix is nonzero. (We also must check that 
all equality constraints are satisfied.) Since $\bv$ is in the Hilbert basis
and the quadrilateral conditions hold,
$S$ is connected by Lemma~\ref{lem302}.
\end{description}
(4b) and (4c) were shown
to take polynomial time in the proof of Theorem~\ref{Nth81}. $~~~\qed$ 

\vspace*{+.1in}
Using Haken's Theorem~\ref{thm302} instead of the Jaco-Tollefson 
result 
we still obtain the weaker result that the {\sc unknotting problem} is in 
the polynomial hierarchy $\Sigma_2^p$ \cite{Stockmeyer:76}:
One uses the same certificate as above,
except for the proof that $S$ is connected. For this step, we use a {\bf co-NP}.
oracle that answers the question: Is $\bv$ in the Hilbert basis of 
$\sC_M$? This problem
is in {\bf co-NP}, because one can certify that $\bv$ is not in the Hilbert 
basis by guessing $\bv_1 ,
\bv_2 \in \sC_M$ such that $\bv = \bv_1 + \bv_2$. This puts the 
problem in the complexity class 
$\mbox{\bf NP}^{\mbox{\bf co-NP}} = \mbox{\bf NP}^{\mbox{\bf NP}}
= \Sigma_2^p$.
A stronger result than this cannot be expected, since it is
known that the problem of testing whether a vector $\bv$ is in the Hilbert 
basis of a pointed
rational cone $\sC$ is {\bf co-NP}
complete, see Seb\"{o} 
\cite[Theorem~5.1]{Seb90}.

\section{Certifying Splittability}\label{csp} 
We treat the {\sc splitting problem} with a modification of the method 
described above.
We use the criterion for splittability of a link given in 
Theorem~\ref{thm504}. The construction of the
certificate, and its verification, take place in the following steps. 

\paragraph{Splitting Link Certificate.}
\begin{enumerate}

\item Given a link diagram ${\cal D}$,
construct a piecewise-linear link $L$ in $\RR^3$ that has regular 
projection ${\cal D}$. From
it construct a good triangulation of $M_L \cong {\bf S}^3 - R_L$ which 
contains $t$ tetrahedra, with $t = O(n)$, and with a meridian marked
in each component of $\partial M_L$. (Use Lemma~\ref{lem501}.) 

\item Guess a suitable vertex solution $\bv \in \ZZ^{7t}$ to the Haken 
normal equations for $M_L$.
(This solution can be written in polynomial length by 
Lemma~\ref{lem401}.) Verify the quadrilateral disjointness conditions. 
Let $S$ denote the associated normal surface, so $\bv = \bv (S)$. 

\item Verify that $S$ is a sphere that splits two components of $\partial 
M_L$. \begin{enumerate}
\item Verify that $S$ is connected by verifying that $\bv$ is a minimal 
vertex solution.

\item Verify that $S$ is a sphere by verifying that $\chi(S)=2$. 

\item Verify that $S$ separates two components $T$ and $T'$ of $\partial 
M_L$ by verifying that the number of intersections of $S$ with the marked 
arc joining $T$ and $T'$ is odd.
\end{enumerate}
\end{enumerate}

We obtain an algorithm for testing a link diagram for splittability by 
exhaustive search for a certificate of the kind above. The {\em Splitting 
Link Algorithm} searches the list $\sL_V$ of minimal vertex solutions 
given by (\ref{eq602}). 

\begin{theorem}\label{thm901}
There is a constant $c''$ such that the Splitting Link Algorithm decides for 
any $n$-crossing link diagram whether it represents a splittable link 
using at most $O( \exp (c'' n))$ time and $O(n^2 \log n)$ space, on a Turing 
machine. \end{theorem}

\paragraph{Proof.}
The correctness of the algorithm follows from the list $\sL_V$ being 
exhausting by Theorem~\ref{thm505} and from the splittability criterion 
Theorem~\ref{thm504}. The connectivity of $S$ and $\chi (S) =2$ imply 
that $\partial S = \emptyset$ and that $S$ is a sphere.

Each part of step~(3) can be verified in polynomial time using 
$O(n^2 \log n)$ space, by the arguments in the proof of 
Theorem~\ref{Nth81}. For step~(3c), exhaustively check marked arcs 
which between them connect every pair of components of the link.
There are at most $O(n^2)$ such marked arcs, each containing at most 
$O(n)$ edges in $M_L$, and we compute the intersection number 
$(\bmod~2)$ of each arc with $S$, as in Theorem~\ref{Nth81}. If $S$ 
certifies a splitting of the link, then at least one arc will have odd 
intersection number with $S$. 

The exponential running time bound follows from the cardinality of 
$\sL_V (M_L)$ in (1) of Lemma~\ref{lem402}.~~~$\qed$ 

\paragraph{Proof of Theorem~\ref{linknp}.} The existence and correctness 
of the Splitting Link Certificate follows from Theorems \ref{thm504} and 
\ref{thm505}. The polynomial-time verifiability follows from the proof of 
Theorem~\ref{thm901} above.~~~$\qed$ 

\section{Determining the Genus}\label{dg} 
The unknotting algorithm of Section 8
can easily be generalized to solve the {\sc genus problem}
in polynomial space. 

\subsection*{Genus Algorithm}
\begin{enumerate}

\item Given a link diagram ${\cal D}$ verify that it is a knot diagram 
$\sK$.

\item As before, construct a piecewise-linear knot $K$ in $\RR^3$ that 
has regular projection ${\cal K}$, together with a good triangulation 
of $M_K \cong {\bf S }^3 - R_K$
which contains $t$ tetrahedra, with $t = O(n)$, and with a meridian 
marked in $\partial M_K$.

\item
Construct an exhausting list $\sL_H \subset \ZZ^{7t}$ that includes all 
fundamental solutions $\bv \in \ZZ^{7t}$ to the Haken
normal equations for $M_K$. This list is taken to be $$\sL_H = \{ \bv = 
(v_1, v_2, \ldots, v_{7t} ) \in \ZZ^{7t}: 0 \le v_i \leq t 2^{7t+2} \} ~.
$$

\item
For each $\bv \in L$ test if $\bv$ is admissible. If so, let $S$ denote a 
normal surface with $\bv = \bv (S)$. The following steps test if $S$ is a 
two-sided connected surface and, if so, compute its genus $g(S) = g( \bv 
)$. 

\item
\begin{enumerate}
\item\label{conna} Verify that $S$ is connected by verifying the 
connectedness of an undirected graph with nodes corresponding to 
triangles in the triangulation of $S$ and edges joining matching triangles. 
Otherwise go to the next $\bv$.

\item\label{connb} Verify that $S$ is orientable by
verifying the non-connectedness of an
undirected graph with nodes representing each of the 
two sides of triangles in the triangulation and edges joining matching 
sides of matching triangles.
(Since the surface $S$ is connected and is an embedded surface in an 
orientable manifold,
$S$ is orientable if and only if it is two-sided.) 

\item\label{connc} Verify that $\partial S$ is non-empty and connected 
(as an undirected graph), so that it is topologically $\bS^1$.
Otherwise go to next $\bv$.

\item Verify that $\partial S$ is a longitude in $\partial M_K$ by 
verifying that the homology class
$[\partial S] = (0,\pm 1)$ in $H_1(\partial M_K; \ZZ /2 \ZZ)$. 
Otherwise go to next $\bv$.

\item
Compute the genus $\displaystyle g(S) = \frac{1-{\chi (S)}}{2}$.
\end{enumerate}

\item
Output the minimal value of $g(S)$ found in step (4). \end{enumerate}

We obtain the following complexity bound for this algorithm. 

\begin{theorem}\label{thm1001}
There is a constant $c'''$ such that for any $n \ge 1$ the Genus 
Algorithm computes the genus $g(K)$ of the knot $K$ represented by a 
given $n$-crossing knot diagram $\sK$ and runs in time at most $O(2^{c ''' 
n^2})$ and uses space at most $O(n^2)$ on a Turing machine. 
\end{theorem}
\paragraph{Proof.}
The correctness of the Genus Algorithm follows from Corollary~ 
\ref{cor502}. It remains to bound the time and space requirements of each 
step of the algorithm. 

In running through the list $\sL_H$, we do not attempt to recognize 
which vectors $\bv$ give fundamental solutions. We simply compute the 
genus for each of them, whenever it is possible. We go through the list 
$\sL_H$ in lexicographic order. The key point is to show each step 
requires polynomial space. In Steps~(\ref{conna}),~(\ref{connb}) and 
(\ref{connc}), we use the fact that in an undirected graph in which nodes 
can be written down in polynomial length and in which adjacency of nodes 
can be tested in polynomial space, the connectedness of the graph can be 
determined in polynomial space (see Savitch \cite{Sav:70}).

In step (5d), we compute the intersection number of $\partial S$ with a 
marked meridian and longitude in $\partial S$. We trace the curve $S$, 
assigning an orientation to each segment, and keep a running total of the 
intersection number. This can be done in polynomial space.

All other steps can be computed in polynomial time, as in 
Theorem~\ref{Nth81}. The resulting time bound is dominated by the 
number of elements $O(2^{c_1 t^2})$ in $\sL_H$, given by the proof of 
Lemma~\ref{lem402}. A bound of $O(2^{c_2 t^2})$ also occurs in the 
connectivity algorithms in (5a) and (5b).
The space bound is dominated by steps (5a) and (5b), and is $O(n^2 )$.~~~$\qed$

\paragraph{Proof of Theorem~\ref{genusps}.} This is 
immediate from Theorem~\ref{thm1001}.~~~$\qed$ 

\section{Conclusion}\label{conc}
We know of no non-trivial lower bounds or hardness results for any of the 
problems we have discussed; in particular, we cannot even refute the 
implausible hypothesis that they can all be solved in logarithmic space. 
There are also a great many other knot properties and invariants apart 
from those considered here, and for many of them it is a challenging open 
problem to find explicit complexity bounds. 

One interesting question is whether the {\sc unknotting problem} is in {\bf 
co-NP}. Thurston's geometrization theorem for Haken manifolds can be
used to show
that knot groups are residually finite \cite{Hempel:92}. It follows that a 
non-trivial knot has a representation into a finite permutation 
group with non-cyclic image.
Unfortunately no way is yet known to bound the size of this group; 
if the number of symbols in the smallest such permutation group were 
bounded by a polynomial in the number of
crossings, then the {\sc unknotting problem} would be in {\bf co-NP}. In 
practice the order of such a group seems to be quite small.

Perhaps the most important of the open problems is to determine the 
complexity of the {\sc knot equivalence problem} 
As mentioned before, a
decision procedure is known (see Waldhausen \cite{Wald:78} and 
Hemion \cite{Hem:92}). However at present it is not even
clear whether the resulting algorithm is primitive recursive. 

\noindent
{\bf Acknowledgments.}
The authors thank J. Birman, R. Pollack, P. Shor, J. Sullivan and G. M. Ziegler for helpful 
comments and references. The second author thanks J. Birman
for introducing him to this problem and to J. Hass during the
MSRI special year on Low-Dimensional Topology.

{\tt email: hass@math.ucdavis.edu} \\
\hspace*{.85in}{\tt jcl@research.att.com} \\ \hspace*{.85in}{\tt 
nicholas@cs.ubc.ca}
\end{document}